\documentclass[11pt, a4paper]{article}
\usepackage{subcaption}
\usepackage{latexsym}
\usepackage{indentfirst}
\usepackage{graphicx}
\usepackage{amsmath}
\usepackage{amssymb}
\usepackage{accents}
\usepackage{tikz}

\newcommand{\ubar}[1]{\underaccent{\bar}{#1}}

\usepackage{hyperref, authblk}
\usepackage[mathlines]{lineno}

\begin{document}
\title{Mathematics for energy systems:\\ Methods, modeling strategies, and simulation}
\author{Nicklas Jävergård, Grigor Nika, Adrian Muntean\thanks{Corresponding author: \url{adrian.muntean@kau.se}}}
\affil{Department of Mathematics and Computer Science \protect\\ Karlstad University, Sweden \protect\\ Universitetsgatan 2\protect\\ 65188 Karlstad, Sweden}

\date{}
\maketitle

\begin{abstract} We offer an insight into our mathematical endeavors, which aim to advance the foundational understanding of energy systems in a broad context, encompassing facets such as charge transport, energy storage, markets, and collective behavior. Our working techniques include a combination of well-posed mathematical models (both deterministic and stochastic), mathematical analysis arguments (mostly concerned with model dimension reduction and averaging,  periodic homogenization), and simulation tools (numerical approximation techniques, computational statistics, high-performance computing). 
\end{abstract}

\textbf{MSC 2020: } 91A16,  65M08, 62R07, 00A71, 35B27

\textbf{Key words}:  Control theory, mean-field games, partial differential equations, statistical aspects of data science, mathematical modeling, homogenization


\section{Introduction}
There is no doubt that many of the societal problems faced in the world today relates to matters of energy, either directly or indirectly. 
Some of these can be addressed by technological advances whilst others deal with human action and behavior. 
Our goal is to help the understanding of these topics broadly seen using mathematical modeling, analysis and simulation.

In this paper we divulge some topics that we collected under the umbrella of mathematics for energy systems. 
We present this collection of topics from an applied analysis perspective. 
The work presented herein spans a rather broad view. 
Some models concern the subject of control and Mean-Field Games (MFGs). MFGs describe a population of infinitesimal, anonymous and indistinguishable agents who interact through their distribution. Players minimize their cost functional which may be constructed to achieve some macroscopic goal for instance, grid stability, minimize energy consumption or maximize self-sufficiency.
Other address multiscale materials science questions, trying to understand the production, transport and storage of energy. 

\begin{figure}[!htb]
    \centering
    \includegraphics[scale=0.45]{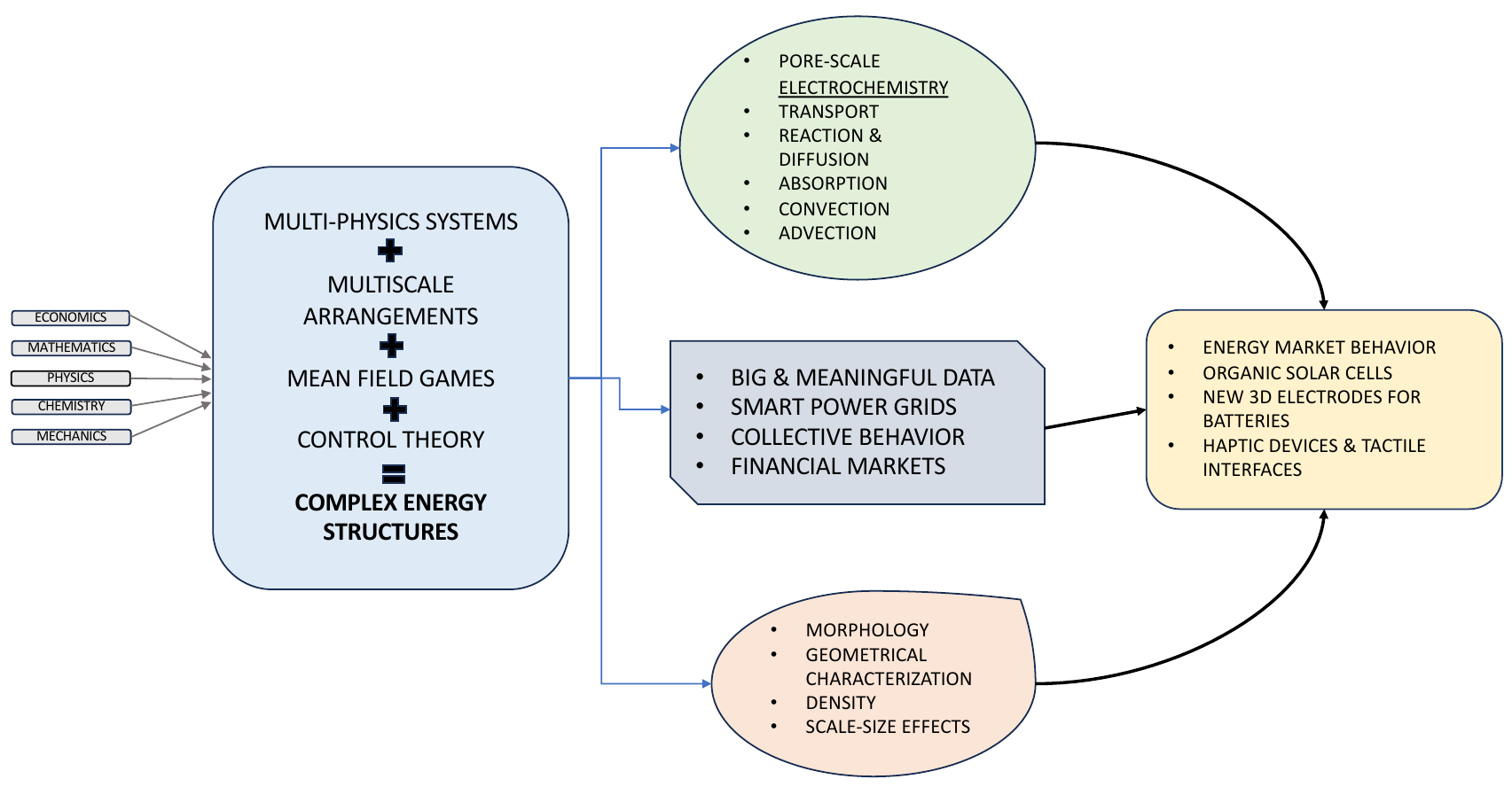}
    \caption{Synergetic mathematical approach to modeling energy systems, including materials for energy harvesting and storage as well as collective behavior in relation to such systems.}
    \label{fig:CorrPlot}
\end{figure}

We present our approach visualized in Figure \ref{fig:CorrPlot}. This is not yet a complete story, many questions still remain to build a realistic description. Currently, we are looking for ways to connect the MFG theory with multi-scale material science to address how units for storage or harvesting should be constructed, see \cite{Welch2019}.

\section{Synthetic data generation}\label{sec:SDG}
This topic rises naturally when trying to validate models in other projects. Specifically, in the case of MFGs (details in section \ref{sec:MFG}), we had need of data detailing the consumption and production of electricity on an aggregate level. Data of this kind exists and would be highly valuable to us as researchers to test models and ideas, however this data is often protected due to issues of privacy or competition. To resolve this we construct a method for generating synthetic data that respect the correlations between features of a dataset without letting the end user of the generated data infer to much about the individual data points in the original dataset. There are numerous ways of generating synthetic data for different purposes. Methods in the field of machine learning for expanding existing or generate synthetic datasets using Generative Adversarial Networks (GANs) and the like, see \cite{Alvaro2022}. There are probabilistic approaches for representing statistical information in datasets as in \cite{Gogoshin2020}. Our method is a probabilistic white box approach, that is we calculate the empirical distribution of each feature. As well as the conditional distributions of any doublet and triplet connection between features of a dataset. This information is then used to draw new representative observations that is stored as a synthetic dataset.

We consider a dataset $\mathcal{O}$, consisting of $N_f$ columns or herein called features and denoted $f_i$, for $i\in\{1, 2, \dots, N_f\}$. Each feature contains $M$ observations, thus $\mathcal{O}=[f_1, f_2, \dots, f_{N_f}]\in \mathbb{R}^{M\times N_f}$. We express the grid points of our discretization in \eqref{eq:discrete}, viz.

\begin{equation}
    x_i^n := \min\{f_i\} + \frac{\max\{f_i\} - \min\{f_i\}}{N}\cdot n \quad \text{ for }n \in \{0, 1, \dots, N\}.
    \label{eq:discrete}
\end{equation}
Using \eqref{eq:discrete}, we define the intervals $A_i^n := [x_i^{n-1}, x_i^n)$ for $n \in\{1,..., N-1\}$ and $A_i^N := [x_i^{N-1}, x_i^N]$ for each feature $f_i$. To facilitate the overall description, we introduce the following convention what concerns the use of the indices: Every time we use them from this point on, we mean that 
\begin{center}
$i, j,k \in\{1, \dots, N_f\}$ and $n, m, \ell\in \{0, \dots, N\}$. 
\end{center}
We are concerned with computing the following probabilities and conditional probabilities, respectively:
\begin{equation}
\begin{aligned}
    p_{i, n} &= P(A_i^n),\\
    p_{[(j,m)|(i, n)]} &= P(A_j^m|A_i^n), \\
    p_{[(k, \ell)|(i,n),(j,m)]} &= P(A_k^\ell|A_i^n, A_j^m),
    \label{eq:teoretical}
\end{aligned}
\end{equation}
for all combinations of indices $(i,j,k)$ and $(n, m, \ell)$ such that $i\neq j \neq k$. These are approximated using the original dataset. Let $f_i^s$, $s\in  \{1, 2, \dots, M\}$ be an entry in feature $f_i$, using this notation we may estimate \eqref{eq:teoretical} as in 

\begin{equation}
\begin{aligned}
    p_{i, n}\approx \hat{p}_{i, n} &= \frac{1}{M}\sum_{s=1}^M\chi_{A_i^n}(f_i^s),  \\
    p_{[(j,m)|(i, n)]}\approx \hat{p}_{[(j,m)|(i, n)]} &= \frac{\sum_{s=1}^M\chi_{A_i^n\times A_j^m}(f_i^s, f_j^s)}{\sum_{s=1}^M\chi_{A_i^n}(f_i^s)}, \\
    p_{[(k, \ell)|(i,n),(j,m)]} \approx \hat{p}_{[(k, \ell)|(i,n),(j,m)]} &= \frac{\sum_{s=1}^M \chi_{A_i^n\times A_j^m \times A_k^\ell}(f_i^s,f_j^s, f_k^s)}{ \sum_{s=1}^{M} \chi_{A_i^n\times A_j^m}(f_i^s, f_j^s)}.
    \label{eq:estimate}
\end{aligned}
\end{equation}
In \eqref{eq:estimate} $\chi_{A_i^n}(f_i^s)$ is 1 if $f_i^s\in A_i^n$ otherwise it is 0. We generate synthetic data using these conditional distributions. More details are given in our preprint \cite{jävergård2024}, where the proposed method is applied to a real dataset concerning energy consumption of households in the Madeira Islands. Here we present a constructed example to give an idea of what to expect.

\begin{figure}[!htb]
    \centering
    \includegraphics[width=12cm]{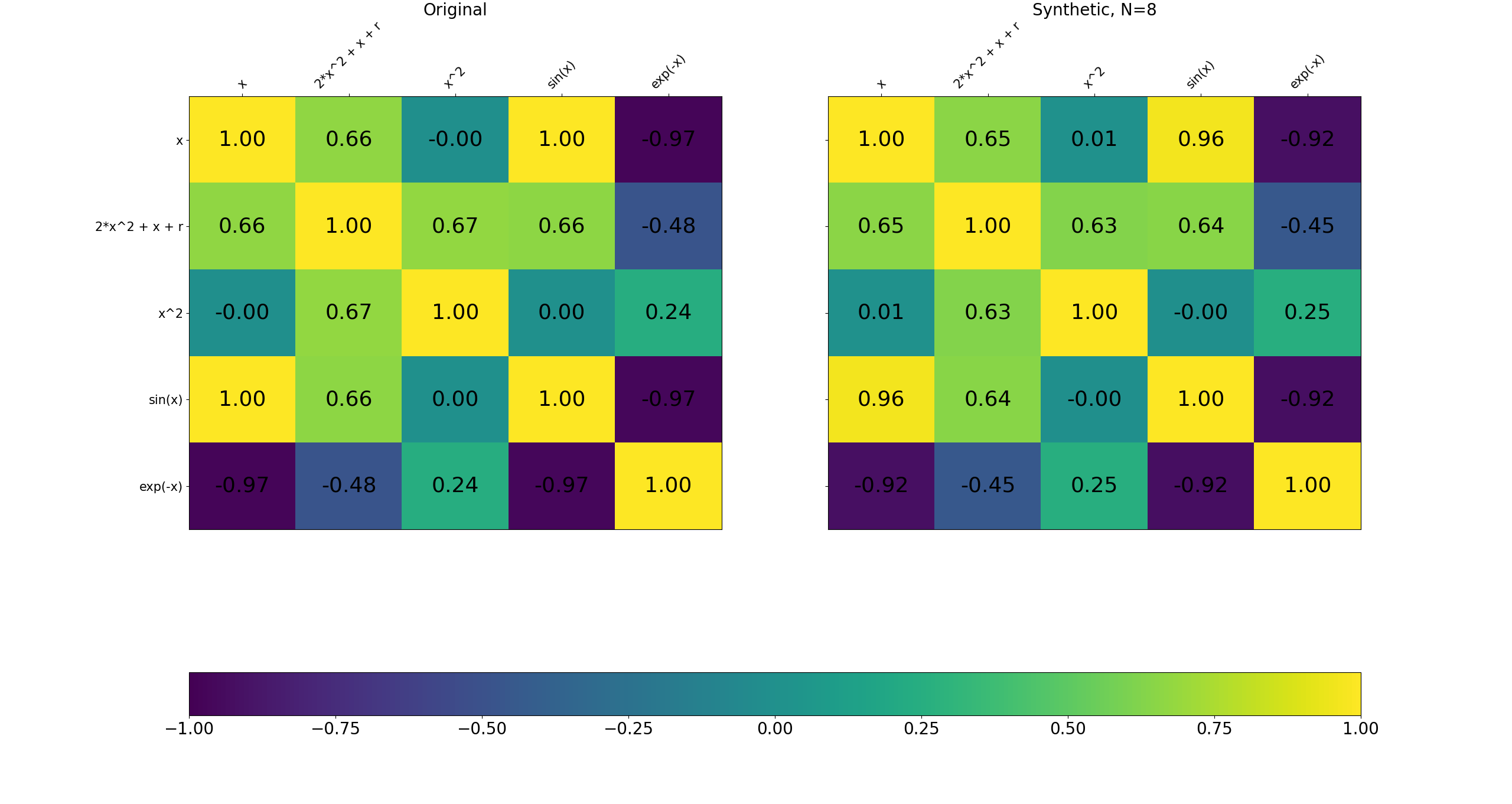}
    \caption{Pearson's correlation coefficient for the original dataset (left) and the synthetic dataset (right) for $N=8$. For information on how this method preserves distributions, see \cite{jävergård2024} for an application to a different dataset.}
    \label{fig:CorrelationPlot}
\end{figure}

\noindent In our manufactured example, we take  $N_f=5$, and $M=1000$. Let $x$ be a vector containing $M$ equidistant values spanning the interval $[-1, 1]$. The features of the dataset are:
\begin{equation}
    \begin{aligned}
        f_1(x) &= x, \\
        f_2(x) &= 2x^2 + x + r, \\
        f_3(x) &= x^2, \\
        f_4(x) &= \sin(x), \\
        f_5(x) &= \exp(-x),
    \end{aligned}   
\end{equation}
where $r$ is a random vector, such that each entry belongs to a uniform distribution on the interval $[0, 1]$. In Figure \ref{fig:CorrelationPlot} we show the Pearson correlation coefficient for the original dataset and the synthetic one for a particular choice of $N$. As can be seen in Figure \ref{fig:CorrelationPlot} we are able to fairly represent the correlations of this simple dataset using the synthetic one. 

The ultimate goal would be to arrive at a point where the privacy and fidelity of the synthetic data can be tuned using the resolution and depths of conditional probabilities to suit the needs of both the data owner (concerns of privacy) and that of the data user (researcher). To this end one could apply differential privacy as a measure as described, for instance, in \cite{Boedihardjo2024}.

\section{MFG for distributed cooling} \label{sec:MFG}
The topic of MFGs was introduced in \cite{Lasry2007}, \cite{Malhame2006} by Lasry, Lions, Malham\'e, Huang and Caines motivated very much by applications in economics.  
The model we briefly describe in this section is a simple MFG model of how one could achieve peak shaving or valley filling on the electrical demands on a large electrical grid using control of thermal devices. We refer the reader to \cite{Bagagiolo2014} for more details and mathematical analysis of this model.
The basic idea is to use the thermal capacity of these devices to de-synchronize the power consumption whilst having approximately the desired temperature for all agents. Pick $\ubar{x}, \Bar{x} \in \mathbb{R}$ such that $\Bar{x}>\ubar{x}$ and define $\Omega := (\ubar{x}, \Bar{x})$. Let $x\in \Omega$ be a real number, representing temperature. Similarly pick $0<T \in \mathbb{R}$ and let $S :=(0, T)$, where $T$ is the terminal time of the model. The model consists of three major parts. The Kolmogorov-Fokker-Planck equation governing the evolution of the distribution of temperatures in the system, The Hamilton-Jacobi-Bellman (HJB) equation whose solution is the value function and the control function which attains the supremum in the HJB equation when optimal. Starting with the Kolmogorow-Fokker-Planck equation, take $m : S \times \Omega \rightarrow \mathbb{R}$ sufficiently smooth such that is solves the following system 
\begin{align}
    &\frac{\partial m}{\partial t} + \frac{\partial}{\partial x}\bigg(f(x, u(x, t))m(x, t)\bigg) = 0  \quad \text{ in } S\times \Omega, \label{eq:M}\\
    &m(\ubar{x},t) = m(\Bar{x}, t) = 0 \quad \text{ for all } t\in S, \\
    &m(x, 0) = m_0(x) \quad \text{for all } x \in \Omega, \\
    & \int_{\ubar{x}}^{\Bar{x}} m(x, t)dx = 1 \quad \text{for all } t \in S,
\end{align}
where $m_0 : \Omega \rightarrow \mathbb{R}$ is a sufficiently smooth initial condition and $f(x, u^*(x, t)) : \Omega \times [0,1] \rightarrow \mathbb{R}$ given by

\begin{equation}
    f(x, u^*(x, t)) := -\alpha x+\sigma u(x,t) + c,
    \label{eq:MFG_velocity}
\end{equation}
for $\alpha>0$, $\sigma = -\alpha(\Bar{x}-\ubar{x})$, and $c=\alpha \Bar{x}$. The control function $u : S\times\Omega \rightarrow [0, 1]$ described the action of a agent and controls the drift in the equation for $m(x, t)$ when the referring to the optimal control we use $u^*$. The optimal control is the control function that minimize the integral over the cost function $g : \Omega^2\times [0, 1]\rightarrow \mathbb{R}$ in \eqref{eq:MFG_cost} using the notation $\Bar{m}^+:= \frac{\Bar{m} + |\Bar{m}|}{2}$ and $\Bar{m}^- := \frac{\Bar{m} - |\Bar{m}|}{2}$.

\begin{equation}
    g(x, u, \Bar{m}) = ru + qx^2 + h\Bar{m}^+u + k \Bar{m}^-(1-u),
    \label{eq:MFG_cost}
\end{equation}
for $r, q, h, k>0$ and the mean value of $m(x, t)$ is $\Bar{m}(t)$. In \ref{eq:MFG_cost} the reference temperature is taken to be 0.
The terms in \eqref{eq:MFG_cost} represent the cost of the device being on, the cost of deviating from the reference temperature, the cost of cooling when the mean temperature is above the mean reference temperature, and the cost of not being on when the mean temperature is low. The last two terms direct the de-synchronization of the players (devices). To minimize the cost of operating within this system as dictated by \eqref{eq:MFG_cost} we need to compute the optimal control as given in \eqref{eq:MFG_control}

\begin{equation}
    u^*(x, t) = \underset{s\in[0,1]}{\text{argmax}}\{-f(x, s)\frac{\partial v}{\partial x}(x, t) - g(x, s, \Bar{m})\} \label{eq:MFG_control}.
\end{equation}

\begin{figure}
    \centering
    \includegraphics[width=0.8\linewidth]{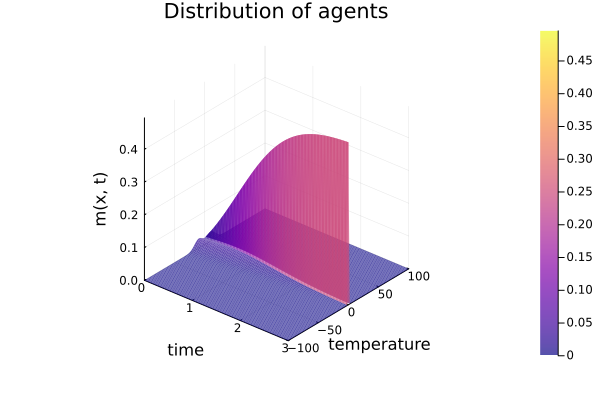}
    \caption{Time evolution of distribution of temperature, $m(x, t).$}
    \label{fig:M}
\end{figure}

\begin{figure}
    \centering
    \includegraphics[width=0.8\linewidth]{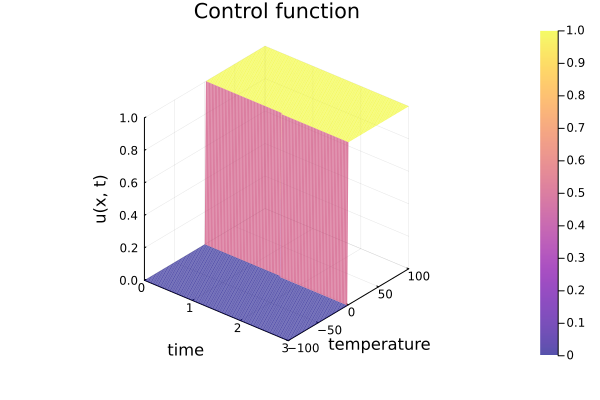}
    \caption{Time evolution of the control function of agents.}
    \label{fig:U}
\end{figure}

Lastly we take a look at the Hamilton-Jacobi-Bellman equation, find $v : S \times \Omega \rightarrow \mathbb{R}$ sufficiently smooth such that
\begin{align}
    -  &\frac{\partial v}{\partial t} + \sup_{u\in U}\{-f(x, u)\frac{\partial v}{\partial x} - g(x, u(x, t), \Bar{m})\} = 0\quad \text{in } S\times \Omega, \label{eq:V}\\
     &v(x, T) = \psi(x) \quad \text{ for all } x \in \Omega, 
\end{align}
which is coupled to the Kolmogorov-Fokker-Planck equation via the optimal control, the former is a forward equation and the later is a backward equation in time. We have constructed an Picard iterative scheme using upwind finite volume discretization, for more information on numerical methods within MFGs see \cite{Mathieuinbook}. For the numerical example shown in Figure \ref{fig:M} and Figure \ref{fig:U}, we have taken as initial datum 
\begin{equation}
    m_0(x) = \frac{1}{\sqrt{2\pi\sigma^2}}\bigg(e^{-\frac{(x+\mu_1)^2}{\sigma^2}} + e^{- \frac{(x-\mu_2)^2}{\sigma^2}}\bigg), x\in \Omega, \label{eq:M_init}
\end{equation}
where $\mu_1=-\mu_2=10$ and $\sigma =7$. 

In Figure \ref{fig:M} we see the distribution of agents evolving according to the game presented in \eqref{eq:M} through \eqref{eq:V} for the initial condition given in \eqref{eq:M_init}. 
Associated with this evolution of agents we also compute the control function. The control function for such a game turns out to be simple as can be seen in Figure \ref{fig:U}.

\newpage

\section{Morphology formation in ternary  interacting mixtures}\label{sec:morphology}

In materials science, there is a lot of focus on building efficient materials to harvest green energy.  In this context, we are mostly interesting in exploring by means of mathematical modeling and simulation tools a couple of questions related to the formation of the internal morphology  used for organic solar cells and link this morphology with their efficiency as device; we refer the reader to  \cite{Moons} for a discussion of the materials science context of polymer photovoltaic cells - the main target here. Essentially, we model in \cite{Andrea_EPJ} the interacting mixture between two polymers and a solvent by means of a stochastically interacting particle system. 

Denote by $\sigma: \Omega \longrightarrow \{-1,0,+1\}$ the spin of a particle at site $x \in \Omega$. In this context, the triplet $(-1,0,+1)$ refers to polymer $A$, solvent $S$, and polymer $B$, respectively. In the simulations, we are going to show in this section the space occupied by these elements will be colored blue, red, and respectively, yellow.   We endow the system  with a Blume-Capel dynamics governed by the Hamiltonian: 
    \begin{equation}\label{Eq_HamAndrea}
        H(\sigma) := \frac{1}{2} \sum_{\substack{    x \neq x' \in \Omega : \\ 
        |x-x'| =1}} M_{\sigma(x),\sigma(x')}, \mbox{ where } M := \begin{bmatrix}
        0 & 1 & 4\\
        1 & 0 & 1\\
        4 & 1 & 0
    \end{bmatrix}
    \end{equation}
  is the so-called interaction matrix. Each entry in the matrix tells how repulsive the corresponding inter-particle interaction is. 
In previous works (cf. e.g. \cite{Andrea_PhysRevE} and references cited therein), Monte Carlo simulations were used to study the formation of morphologies and their growth. The main scientific challenge is to unveil the effect of the evaporation of one of these species (the solvent) on the growth of the formed phases. 
\begin{figure}[h]
    \centering
    \includegraphics[width=0.45\linewidth]{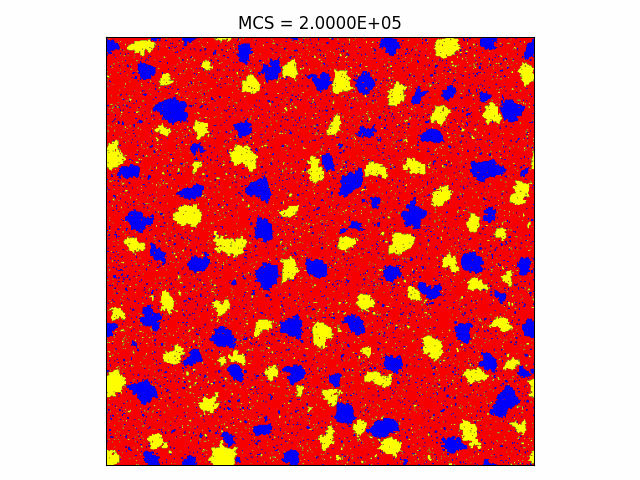}
    \caption{Typical Monte Carlo simulation output for a solvent level of $0.8$.}
    \label{fig:MC}
\end{figure}
It turns out to be conceptually very difficult to do such simulations in $2D$ while taking evaporation into account.  The $3D$ case is even more complex. To avoid such obstacles\footnote{The hinder is both from the modeling and computational points of view. What concerns the modeling part, it is not clear to us cut how one can implement naturally even a simple evaporation dynamics in a (conservative) Monte Carlo setup.}, we became interested in the continuum (hydrodynamic) limit for such interacting particle systems (very much in the spirit of the monograph \cite{Presutti}).  Using the results reported in \cite{Marra}, we studied in 
\cite{lyons2024phase, lyons2023continuum} a continuum counterpart of an interacting particle system close to what is captured by the Hamiltonian \eqref{Eq_HamAndrea}. It is about the following coupled parabolic system with nonlinear and nonlocal drifts. The choice of nonlinearity in the structure of the drift arises from a particular regularized Hamiltonian (a simplification of \eqref{Eq_HamAndrea} as explained in \cite{Marra}). Essentially, we are looking for a pair of functions $(m,\phi)$ that satisfy the following system of evolution equations posed in the domain $\Omega$:
  \begin{gather}
    \left\lbrace
    \begin{aligned}
        \partial_t m &= \nabla \cdot \left[\nabla m - 2 \beta (\phi -m^2 ) (\nabla J * m) \right] \mbox{ in } (0,T)\times\Omega\\
        \partial_t \phi &= \nabla \cdot \left[ \nabla \phi - 2 \beta m (1 - \phi) (\nabla J * m) \right]\mbox{ in } (0,T)\times\Omega.
    \end{aligned}\right. \label{SYS}
    \end{gather}
    Both unknowns are requested to satisfy periodic boundary conditions and a suitable initial condition is imposed. Here $\beta$ is the reciprocal of the  temperature, while $J(\cdot)$ is a smooth convolution kernel. It is straightforward to see that if one picks in \eqref{SYS} $\phi=1$ everywhere in $(0,T)\times \Omega$, then the original system reduces to a scalar, nonlinear, and  nonlocal Cahn-Hilliard-type equation. If we think of the morphology formation problem for organic solar cells, then $1-\phi$ means the solvent fraction, while the magnetization $m$ describes the mean spin. Together with $\phi$, $m$ tells what constituents are present at certain position at a given time. 

\begin{figure}[h]
    \centering
    \includegraphics[width=0.45\linewidth]{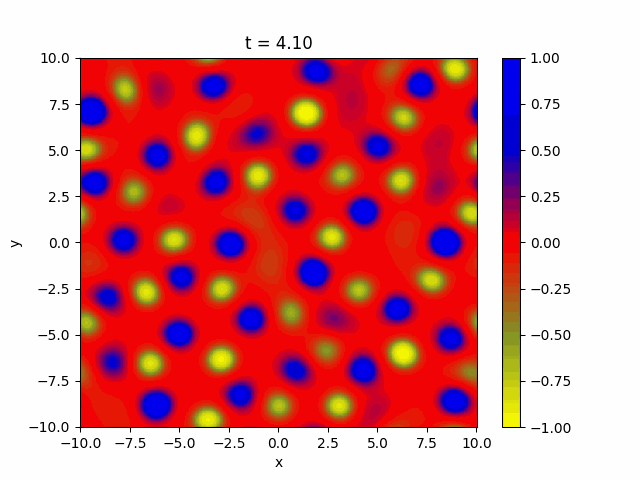}
    \caption{Typical PDE-based simulation output for a solvent level of $0.8$.}
    \label{fig:PDE}
\end{figure}

A typical output of both models is shown Figure \ref{fig:MC} (for the Hamiltonian \eqref{Eq_HamAndrea} with the Kawasaki dynamics) and in Figure \ref{fig:PDE} (for a finite-volume approximation of the solution to the system \eqref{SYS}).

The details on the numerical schemes used to produce these Figures are reported in \cite{Andrea_EPJ} and \cite{lyons2023continuum}, respectively. See also  \cite{lyons2024bound} for related results.

\section{Charge transport through heterogeneous media}\label{sec:charge}
Rising demand for mobile power sources ranging from electric vehicles to portable electronics has spurred substantial research aimed at creating advanced, high-performance electrochemical energy storage devices (\cite{ray2012rigorous}, \cite{vo2019corrector}), large area organic light emitting diodes and organic semiconductors (\cite{glitzky2021existence}, \cite{glitzky2022analysis}), electro-magneto rheological fluids (\cite{vernescu2002multiscale}, \cite{NV20}), and many more. As a result, advanced modeling and simulation techniques become indispensable in understanding, e.g., the transport mechanism of ions in lithium ion batteries.

The mathematical modeling framework for these systems comprises three principal classes of equations: 1) The Nernst-Planck equation, which characterizes the flux and transport dynamics of charges, ions, or electrons within the system. 2) Maxwell’s equations, which dictate the behavior of static electric and magnetic fields\footnote{In a quasi-static regime Maxwell's equations with negligible magnetic fields reduce to an electric potential that satisfies Poisson's equation. Starting from Maxwell's equations for the electrical induction and electrical field, respectively
\begin{align*}
\nabla \cdot D &= \rho q \text{ in }\mathbb{R}^3, \\
\nabla \times E &= 0 \text{ in }\mathbb{R}^3.
\end{align*} 
With the constitutive law $D = \epsilon_r\epsilon_0 E$, it is derived by noting that the electric field $E$ satisfies $\nabla \times E = 0$ which implies that $E$ may be written as $E = -\nabla \Phi$, for some scalar field $\Phi$ and thus $\nabla \cdot (-\epsilon_r\epsilon_0 \nabla \Phi) = \rho q$, where $\rho q$ is the charge density.}. 3) Constitutive equations that define the properties of the medium in which the charges and/or active fields are embedded.

A paradigmatic model of microscopic electrochemical interactions in an incompressible Newtonian fluid within a domain $\Omega$ with boundary $\Gamma$\footnote{The boundary $\Gamma$ is decomposed as $\Gamma := \Gamma_D \cup \Gamma_N, \quad \Gamma_D \cap \Gamma_N = \emptyset$ to reflect Dirichlet and Neumann boundary conditions, respectively.} is given by the following system of balance equations formulated for the transport of charges. The bulk model equations are:
\begin{gather}
    \begin{aligned}
        \partial_t c^{\pm} + \nabla {\cdot} \Big( \mathbf{v}c^\pm - &D^\pm \nabla c^\pm\\
        - &\frac{D^\pm z^\pm e}{kT} c^\pm \nabla \Phi \Big) &{=}& R^\pm (c^+, c^-) \text{ in } (0,T) \times \Omega, \nonumber
    \end{aligned} \\
    \begin{aligned}
       -\Delta \Phi &{=} \frac{e}{\epsilon_0 \epsilon_r}(z^+c^+ - z^-c^-) \text{ in } (0,T) \times \Omega, \\
       -\eta \Delta \mathbf{v} + \frac{1}{\rho} \nabla p &{=} -\frac{e}{\rho} (z^+c^+ - z^-c^-) \nabla \Phi \text{ in } (0,T) \times \Omega, \\
        \nabla \cdot \mathbf{v} &{=} 0 \text{ in } (0,T) \times \Omega, \\
    \end{aligned}
\end{gather}
The system is completed with the following 
boundary and initial conditions:
\begin{gather}
    \begin{aligned}
        \left( -\mathbf{v}c^\pm + D^\pm \nabla c^\pm + \frac{D^\pm z^\pm e}{kT}c^\pm \nabla \Phi \right) \cdot \nu &{=} 0 \text{ on } (0,T) \times \Gamma, \\
        \nabla \Phi \cdot \nu &{=} \sigma \text{ on } (0,T) \times \Gamma_N, \\
        \Phi &{=} \Phi_D \text{ on } (0,T) \times \Gamma_D, \\
         \mathbf{v} &{=} \mathbf{0} \text{ on } (0,T) \times \Gamma, \\
        c^\pm &{=} c^{\pm,0} \text{ in } \{ t = 0\} \times \Omega.
    \end{aligned}
\end{gather}
Here, $c^\pm$ is the (respective) charge density, $(\mathbf{v}, p)$ are the fluid velocity and pressure, $\Phi$ is the electrostatic potential, $T$ is the absolute temperature, $\eta$ is the kinematic viscosity of the fluid, $\rho$ is the fluid density, $D^\pm$ is the diffusivity (of the respective charge density $c^\pm$), $\sigma$ the surface charge, $z$ the charge number, $e$ the elementary charge, $\epsilon_0 \epsilon_r$ is the di-electrostatic permittivity-relative permittivity, and $k$ is the Boltzmann constant. We denote by $\nu$ the outer unit normal to interfaces. If $\Omega$ is a periodic domain (mimicking the geometry of a regular porous material), averaging techniques offer a natural approach to approximate complex problems and derive effective equations that incorporate both the volume fraction and the morphology of the microstructure (e.g., \cite{allaire2010homogenization, allaire2011erratum}, \cite{ray2012rigorous}). Furthermore, additional correction terms can be computed to refine these approximations using the periodic cell problem (see, e.g., \cite{vo2019corrector}). 

These type of models have well founded applications on emerging fields, e.g., design and fabrication of ${\rm 3D}$ thick electrode structures (e.g., \cite{saleh20183d, sun2019hierarchical, liang2019composite}), magnetic drug targeting (e.g., \cite{grillone2017magnetic}), amelioration of endurance of organic electronics (e.g., \cite{kirch2020experimental, Fuhrmann2020}). Furthermore, when dealing with microstructured materials, two essential questions arise regarding the influence of microstructure on physical properties. First, {\it how does the presence of non-centrosymmetric geometries affect phenomena such as charge or heat transport?} Second, {\it how should we account for effects that arise at the scale of the microstructure itself?} 

Enriched models\footnote{To develop multiscale and multiphysics models that incorporate features like the ones mentioned here and, perhaps, enriched boundary conditions, it is essential to establish a precise theoretical framework grounded in the laws of thermodynamics posed for generalized continua. Initially, constitutive laws should be formulated using a generalized Coleman–Noll procedure. In a subsequent stage, where feasible, these laws should be further refined through {\it forward} homogenization theory, carefully considering morphology, volume fraction, and, notably, the size of the representative volume element. This approach would extend traditional averaging techniques, which typically only address morphology and volume fraction, by explicitly including additional microstructural length scales.} to address the above type of questions (as well as others in this direction) have been known for a while in solid mechanics (see, e.g., \cite{Toupin62, Toupin64}, \cite{MT62, ME68}) with similar type of generalizations in fluids mechanics as well, e.g.,  \cite{fried2006tractions}. Moreover, homogenization techniques for such media have been carried out in a large number of works, e.g., \cite{Forest01}, \cite{nika2024derivation}, \cite{RR22}, \cite{NiMu24}.  However, many relevant questions still remain open especially what concerns model computability and validation against experiment (including e.g. multiscale parameter identification and related inverse questions, distributed control, learning of microstructure information via neuronal networks, and so on).

\section{Comment to the reader}
This note should be seen as outreach. It is an invitation to building new collaborations in the field of applied/applicable analysis [combining modeling, mathematical analysis, numerical simulation and data], driven by scientific questions inspired by successful engineering case studies in materials science, in particular, and in nowadays top engineering technology, in general. 

\section*{Acknowledgments} 
The authors are involved in the Swedish Energy Agency’s project Solar
Electricity Research Centre (SOLVE) with grant number 52693-1. GN and AM gratefully acknowledge financial support from the Knowledge Foundation (project nr. KK 2020-0152).
NAISS (project nr. 2023/22-1283) is acknowledged for providing the computational resources needed for  generating synthetic data.  We thank our collaborators P. Huang (Dalarna) and J. Forsman (CGI) for fruitful initial discussions on the synthetic data topic, and  S.-A. Muntean (Karlstad), E.N.M. Cirillo (Rome), R. Lyons (Boulder), and V. Kronberg (Eindhoven) for inspiration around modeling interacting ternary mixtures (either discrete or continuum).

\small
\bibliographystyle{plain}

\bibliography{ref.bib}

\end{document}